\documentclass[12pt,leqno]{article}
\begin{document}
\title{METHOD OF HIDDEN PARAMETERS AND
PELL'S EQUATION}
\author{S. N. Arteha}
\date{Space Research Institute, Profsoyuznaya 84/ 32,\\
Moscow 117810, Russia. ~ E-mail: sergey.arteha@mtu-net.ru}
\maketitle
\begin{abstract}
Using the representation of numbers, the methods of decreasing the number of 
calculation steps for Pell's equation are developed. The formulae 
relating the natural solution $x_0$ of Pell's equation for $A_0$ with the infinite 
number of minimal natural solutions $x_i$ for other $A_i$ are obtained for various 
cases. The parametric representation of Pell's equation solutions 
are obtained with the help of new "method of hidden parameters".
\end{abstract}
\section{INTRODUCTION}

Pell's equation for natural numbers $A, x, y$ ~(actually 
put forward by Fermat in 1657) 
\begin{equation}
y^2 - Ax^2 = 1 
\end{equation}
has been thoroughly studied (see [1-6,9] and references herein). In the 
subsequent text, we shall consider the minimal nontrivial natural solutions 
$(x,y)\ne (0,1)$ only (each solution is unique; the value of $x$ determines 
the value of $y$). 
To solve (1) for the given natural number $A$ which is not a square, one can use 
several methods: the sequential differences (in essence, Euclidian algorithm: measurement of a 
greater number by a smaller one), the method of Wallis-Brouncker, indian 
cyclic method, the method of continued fraction. We shall not 
consider the method of Wallis-Brouncker and indian cyclic method, since they 
are rather unwieldy for using "a pen and 
a paper only", and they use some exhaustive search: the number of algorithmic 
steps is great enough (the cyclic method can be connected with the method of 
sequential differences and the method of continued fraction: see 
interpretation in [2], Sec.~8.2).

We remind that the sequential differences 
$\Delta_i$ of natural numbers $F_i$ and $f_i$ are 
$$
\Delta_i = F_i - f_i~ ,~~~ F_{i+1} = max{[f_i, \Delta_i]}~ , ~~~ f_{i+1} = min{[f_i, \Delta_i]}~ ,
$$
The determinant of the binary quadratic form
\begin{equation}
b_iY_i^2 - c_iX_i^2 + 2a_iX_iY_i = 1
\end{equation}
is invariant:
\begin{equation}
b_ic_i + a_i^2 = A .
\end{equation}
For this method (which is algorithmic) one can obtain the expression ~
$r_i=b_i-c_i+2a_i$ ~
from (2) at the step $i$. For $r_i>1$ the substitutions are: $X_i=Y_i+X_{i+1},~
Y_i=Y_{i+1}$;~ for $r_i < 0$ the substitutions are: $Y_i=X_i+Y_{i+1}~,~ X_i=X_{i+1}$ 
(the "ultimate form" becomes an identity for $X_n=1,~ Y_n=1$). Then, 
substitutions in the reverse sequence give a solution of (1). 

We remind the method of continued fraction [8]. 
Let $A>0$ be any radicand with $\sqrt{A}=<d_0;d_1,d_2,...,d_L>$ and 
$L=L(\sqrt{A})$, where $\sqrt{A}=<d_0;d_1,...,d_L>$ is the continued fraction 
expansion of $\sqrt{A}$ with period length $L$. If $L$ is even, then all 
positive solutions of (1) are given by $x=B_{Ls-1}, s\ge 1$. If $L$ is odd, 
then all positive solutions of (1) are given by $x=B_{2Ls-1}, s\ge 1$. The 
sequence $\{B_n\}$ is the well known recursive sequence given by $B_{-2}=1,
B_{-1}=0, B_n=d_nB_{n-1}+B_{n-2}$ for $n\ge 0$.

The following identity will be used also: 
\begin{equation}
(y^2 - Ax^2)(Y^2 - AX^2) = (yY \pm AxX)^2 - A(yX \pm xY)^2 .
\end{equation}
We remind the well-known statements [2] which prove the procedure of division 
(for the example of $a^2+b^2$):

(a) the product of two numbers, each of which is a sum of two squares,
is a sum of two squares as well (see (4));

(b) if the number, which is a sum of two squares, is divisible by a
prime which is a sum of two squares, then the
quotient is also a sum of two squares;

(c) if the number, which can be written as a sum of two squares, is
divisible by a number which is not a sum of two squares, then the
quotient has a factor which is not a sum of two squares;

(d) if $a$ and $b$ are relatively prime, then every factor of
$a^2+b^2$ is a sum of two squares.

Pell's equation may seem to be fully studied and no new phenomena can be found.
However, observing the minimal nontrivial solutions of Eq.(1), we see very
great differences in values of $x_i$ with $A_i$. The answer to the following
question is unknown: what cases of $A_i$ are the most difficult ones (with
maximal values of $x_i$)? A maximum of solutions of Eq.(1) for $A$ between
any two nearest squares will be named the local maximum, i.e. for
$k^2 < A_k^i < (k+1)^2$ and $i={\overline{1,2k}}$ we have:
$x^{loc}_{max}(k) > x(A_k^i)$. The absolute maxima are the such solutions of
Eq.(1), that we have $x_{max}^{abs}(A'_n) > x(A_i)$ for all $A_i < A'_n$.
Observations demonstrate, that all local maxima are the cases of some prime
$A = p_i$ and of some double prime $A = 2p_j$, and all absolute maxima (for
$A > 46$) are the cases of prime $A = 4n+1$ (the record-holders are tested
from Internet for $A \le 971853031560$). Thus, the study of particular cases of prime $A$ is
interesting. Since factorization of the double prime corresponding to
local maxima does not present in the full cycle, these double prime can be
named "quasiprime" (we do not have the form (2) with $a_i = 0$, but have
one of the characteristic forms for prime numbers of $A$). Some cases of complex
$A$ correspond to one of the characteristic forms of prime $A$ (i.e. the
case of prime $A$ is rather general one). It is interesting to find methods
for lowering of volume of calculations (in comparison with well knowm
methods) for Pell's equation (1) or for a similar equation with arbitrary right side
($\ne 1$), and to obtain explicit formulas.

We shall use a parametric representation of minimal natural solutions. The 
parameters $l$ and $m$ will denote the basic ones. Knowing the values of $l$ 
and $m$ (positive or zero), we can obtain the minimal positive solution $x_0$ 
for a given number $A_0$ (see Propositions 1,4,7, for example). 
Note that the values of parameters $l$ and $m$ are considerably smaller than 
the value of the solution $x_0$, i.e. we operate with comparatively small 
numbers. To find these parameters, it is suffice to use a smaller number of 
(initial only) steps than that in (full cycle) the method of sequential 
differences, or, than that in the continued fraction method (for example, 
see Section 2). The parameter 
$i$ will denote an additional one ($i=\pm 1,\pm 2,...$). This parameter 
determines some new natural solution $x_i$ (linearly proportional to $i$) for 
some new number $A_i$. Relations (linear in $i$ and $x_0$), which are 
determined by the parameter $i$ for fixed parameters $l$ and $m$, will be 
named "horizontal" relations (see Propositions 2, 5, 8, for example). 
"Vertical" relations are those ones between solutions, for which at least 
one of basic parameters ($l$ or $m$) is not fixed (these relations are 
nonlinear in $l$ and $m$). The method of hidden parameters is introduced to 
find "vertical" relations (in the general case for each special case). The 
method (to find the parameters $l$ and $m$ and solutions $x_i$ for $A_i$) 
consists in the following: first, we insert (instead of initial parameters) an 
additional (greater) number of parameters; then, resolving some system of 
equations, we obtain a solution (in terms of new parameters); finally, 
reducing the problem to the initial number of parameters, the minimal natural 
solutions $x_i$ can be obtained (as an identity) for $A_i$. The method is 
demonstrated in detail to derive the Proposition 3 (and is used for 
Propositions 6, 9). The term "hidden parameters" is suggested, since obtaining 
of the solution (in this procedure) seemed to be impossible at a glance (since 
we seemingly rename parameters only).
 It is obvious, that the relation
between a solution $(x_0, A_0)$ and other solutions $(x_i, A_i)$ does not
connected with composition of forms in cycles (since the determinants are
different for these cases). The general solutions demonstrate the
relationship between Pell's equation (1) and the indeterminate linear
equation of type $ax - by = 1$.

The purposes of this article are as follows:

- the method for considerable decreasing the number of calculation steps is 
demonstrated for Pell's equation (for example, each complicated Fermat case can 
be solved with the help of "a pen and a paper" within a few minutes); 

- by obtaining the solution $x_0$ for $A_0$ one can find the infinite number of 
minimal natural solutions $x_i$ for $A_i$ (for filling out the table of the 
minimal natural solutions ($A_i, x_i$)); 

- the primary aim of this work is to establish some identities for 
Pell's equation in an explicit form. Substituting specific values of 
parameters $l, m, i$ into these identities, the solution ($A_i, x_i$) of (1) 
can always be obtained. 

In Section 2 the case of prime numbers $A=4N+1$ is investigated.
It presents the method of decreasing calculations.
The infinite number of "horizontal" relations in Pell's
equation are obtained for $A_i=a_i^2+b_i^2$. The minimality of solutions
is proved. The basic ideas, how to obtain the formulae of "vertical" 
relations in the case of $A=a^2+b^2$, and the method of hidden parameters 
are considered in detail. In Section 3 the case of prime numbers $A=8N+3$ 
(and the method of decreasing calculations) is considered . The 
connections between $(A_0,x_0)$ and $(A_i,x_i)$ are found. Parametric 
solutions are derived. The case of prime numbers $A=8N+7$ is studied in 
Section 4. In Section 5 the other cases of Pell's equation are discussed. 
Some connections between solutions are demonstrated for 
each of these cases. The algorithm of "inverse calculations" is proposed for 
the sequential differences.
In Section 6 the equation $y^2-Ax^2=-3$ is considered. Section 7 contains 
conclusions. 

\section{THE CASE OF PRIME $A = 4N + 1$}

The case of such prime is of great interest.
Notice, that for all special cases, which were proposed by Fermat (namely
$A=61; 149; 109; 433$), the prime numbers are $4N+1$. Besides, conclusions 
from this Section can be applied to the case of some complex $A=a^2+b^2$. \\
~\\
{\Large Proposition 1}. {\it Let} $A$ {\it be a prime number} $A=4N+1=a^2+b^2$ 
($b$ {\it denotes an odd number), then the minimal natural solution of Pell's 
equation (1) is}
\begin{equation}
x = 2|2blm + a(l^2 - m^2)|(l^2 + m^2) ,
\end{equation}
{\it where} $l$ {\it and} $m$ {\it are taken from the "distinctive form"}
\begin{equation}
b(l^2 - m^2) - 2alm = \pm 1 .
\end{equation}

{\it Proof}. The solution of Eq.(1) for a prime number $A=4N+1$ is an even 
number (see congruence modulo 4)
\begin{equation}
x = 2^TSQ
\end{equation}
with coprime odd numbers $S, Q$. Since prime $A$ cannot be
factorized, two possibilities follow: 
$$
(i) ~~~~~~~~~~~~~~~~~~~~~ y - 1 = 2^{2T-j}AS^2 ~, ~~~ y + 1 = 2^jQ^2 ~,
$$
where $j=1$ or $j=2T-1$. Subtracting the first equation from
the second one, we obtain that $2^{j-1}Q^2=2^{2T-j-1}AS^2+1$. 
This result contradicts the minimal nontrivial natural solution (7).
$$
(ii) ~~~~~~~~~~~~~~~~~~~~~ y - 1 = 2^jQ^2 ~, ~~~ y + 1 = 2^{2T-j}AS^2 ~.
$$
Subtracting the first equation from the second one, it follows, that
$2^{j-1}Q^2=2^{2T-j-1}AS^2-1$. This equation is impossible to solve modulo 4 for 
$j=1$. Therefore, as the result ($j=2T-1$), we obtain
\begin{equation}
(2^{T-1}Q)^2 + 1 = AS^2 .
\end{equation}
Since the number in the left-hand side of Eq.(8) has the $a^2+b^2$ type,
the number in the right-hand side must be of identical type (there exists the 
procedure of division). The value $A$ can be represented as $a^2+b^2$. Since 
one of squares in the left-hand side of (8) is unity, the number 
$S$ is a number from primitive Pythagorean triples, that is, $S=l^2+m^2$, 
where $l$ and $m$ are coprime numbers (one of which is odd and one is even). 
The primitive
Pythagorean triple is $S^2=(l^2-m^2)^2+(2lm)^2$. Substituting this into (8), 
we have
$$
(a^2 + b^2)[(l^2 - m^2)^2 + (2lm)^2] = (2^{T-1}Q)^2 + 1 .
$$
Using the identity (4), we obtain
$$
b(l^2 - m^2) - 2alm = \pm 1 ~, ~~~ 2blm + a(l^2 - m^2) = 2^{T-1}Q .
$$
The net result follows from expression (7).
\vskip 0.25in

Using the expression $y = 2A(l^2+m^2)^2 - 1$ and the inequality $l > m$, we see, that
$l \le \sqrt[4]{(y+1)/(2A)} < \sqrt[4]{(\sqrt{A}x+2)/(2A)}$. Therefore, the parameters $l$ and $m$ are comparatively small ones.
We shall demonstrate considerable lowering of calculations for the case of 
$A=61=6^2+5^2$. For comparison, we consider two methods: (I) the method of 
sequential differences; (II) the continued fraction.

(I) The distinctive form (6) is $5(l^2-m^2)-12lm=\pm 1.$ The substitutions are:\\ 
1) $l=m+C \Rightarrow 5C^2-12m^2-2mC=\pm 1;$ 2) $C=m+R \Rightarrow
5R^2-9m^2+8mR=\pm 1;$\\
3) $m=R+D \Rightarrow 4R^2-9D^2-10RD=\pm 1;$ 4) $R=D+V
\Rightarrow 4V^2-15D^2-2DV=\pm 1;$\\
5) $V=D+E \Rightarrow 4E^2-13D^2+5DE=\pm 1;$
6) $E=D+Y_n \Rightarrow 4Y_n^2-3D^2+14DY_n=\pm 1;$\\
7) $D=Y_n+G
\Rightarrow 15Y_n^2-3G^2+8Y_nG=\pm 1;$ 8) $G=K+Y_n \Rightarrow
20Y_n^2-3K^2+2KY_n=\pm 1;$ 9) $K=H+Y_n \Rightarrow 19Y_n^2-3H^2-4HY_n=\pm 1;$
10) $H=Y_n+X_n \Rightarrow 12Y_n^2-3X_n^2-10X_nY_n=\pm 1$.
Therefore, we have ~ $X_n=Y_n=1, H=2, K=3,..., m=21,..., l=58$. Note that 
the full cycle for $A=61$ contains 64 steps (and the appropriate number of 
reverse substitutions). Substituting
parameters $m$ and $l$ into Eq.(5), one obtains the
well known solution $x=226153980$ (instead of usual 64 steps we made only 
10). In a similar manner, the
case of $A=149$ can be solved in 12 steps, 18 steps are necessary for
$A=109$, and so on.
The question of representation of prime $A=4N+1$ as $a^2+b^2$ can be
bypassed (see Section 5 below: for "inverse calculations" we can start from 
the form $3y^2-12x^2-10xy=1$).

(II) Now we consider the method of continued fraction. The continued fraction 
expansion of $\sqrt{A}$ gives the following expression at the step $i$:
\begin{equation}
d_{i-1} + {r_i\over \sqrt{A} + a_i^t}, ~~~ 1\le i\le L+1 . 
\end{equation}
Note that these coefficients determine some form (2). In this case, we have:
$$
a_i = (-1)^{i+1}a_i^t, ~~ b_{2i+1} = b_{2i} = r_{2i}, ~~ c_{2i+2} = c_{2i+1} 
= r_{2i+1}.
$$
The method for lowering of calculations consists in the following: we seek two 
expressions in the continued fraction expansion with $r_j=r_{j-1}$. Therefore, 
we have $a=a_j^t, ~b=r_j$. The parameters $l$ and $m$ is $l=B_{j-1}, ~~m=
B_{j-2}$, where $\{ B_n\}$ is the recursive sequence given by $B_{-2}=1, B_{-1}
=0, B_n=d_nB_{n-1}+B_{n-2}$. Thus,
$$
r_j = r_{j-1} ~ ~ \Rightarrow ~ ~ a = a_j^t, ~~~ b = r_j, ~~~ l = B_{j-1}, ~~~ m = B_{j-2}.
$$
We demonstrate the case of $A=61$.
$$
\sqrt{61}=7+(\sqrt{61}-7)=7+{12\over \sqrt{61}+7}, ~~
{\sqrt{61}+7\over 12}=1+{\sqrt{61}-5\over 12}=1+{3\over \sqrt{61}+5},
$$
$$
{\sqrt{61}+5\over 3}=4+{\sqrt{61}-7\over 3}=4+{4\over \sqrt{61}+7}, ~~ 
{\sqrt{61}+7\over 4}=3+{\sqrt{61}-5\over 4}=3+{9\over \sqrt{61}+5},
$$
$$
{\sqrt{61}+5\over 9}=1+{\sqrt{61}-4\over 9}=1+{\underline{5}\over \sqrt{61}+4}, ~~ 
{\sqrt{61}+4\over 5}=2+{\sqrt{61}-6\over 5}=2+{\underline{5}\over \sqrt{61}+6},
$$
We see that $b=5, a=6, B_0=1, B_1=1, B_2=5, B_3=16, m=B_4=21, l=B_5=58$. 
Substituting $l,m,a,b$ into Eq.(5), the minimal positive solution $x$ follows. 
Instead of 12 steps for the continued fraction expansion of $\sqrt{61}$ we use 
6 steps only, and instead of calculation of $B_{21}$ (it is the solution - see 
Introduction), we calculate $B_5$ only. Therefore, the number of calculation 
steps can be considerably decreased.\\
~\\
{\Large Proposition 2}. {\it Let} $x_0$ {\it be the minimal natural solution of 
Eq.(1) for} $A_0=a_0^2+b_0^2$ {\it and definite parameters} $l, m$ {\it from Eq.(6); 
then} $x_i$ {\it are the minimal natural solutions of Pell's equation for} $A_i$ {\it 
and} $i\ne -1$, {\it where}
\begin{equation}
x_i = |~x_0 + 2iS^3~| , ~~~~~~S = l^2 + m^2, 
\end{equation}
\begin{equation}
A_i = A_0 + i{x_0\over S} + i^2S^2. 
\end{equation}

{\it Proof}. In fact, parameters $l$ and $m$ are fixed ("horizontal" 
relations). Introducing substitutions 
\begin{equation}
k = l^2 - m^2  > 0, ~~~ t = 2lm , 
\end{equation}
we rewrite (6) as
$$
bk - at = \pm 1 . ~ \leqno{(6')}
$$
It follows from this: 
\begin{equation}
b_i = b_0 + it , ~~~ a_i = a_0 + ik , ~~~~~ i=\pm 1, \pm 2,...
\end{equation}
Substituting (13) into expression (5), we obtain $x_i=2|b_it+a_ik|(l^2+m^2)$. Taking
into account substitution (12) and the identity $k^2+t^2=S^2$, one gets the net formulae 
(10) and (11). The  values  ($l$,$m$) are determined by ($k$,$t$). It is 
well known that there exists one minimal positive solution ($k_0,t_0$) of 
Eq. (6') with $|k_0| < |a|, |t_0| < |b|$ for the given $a$ and $b$. Therefore, the same 
values $(k_0,t_0)$ are minimal ones for $a_i$ and $b_i$, $i\ne -1$ (see 
expression (13)). The solution minimality is proved. 
\vskip 0.25in

Note, that if $x_0$ is some nonprimitive solution, then $x_{-1}$ can be 
a primitive solution (if $|a_{-1}| > |a_0|$ and $|b_{-1}| > |b_0|$). For example, we 
can find from the case of $A_0=2$ (parameters $l=2,m=1$ give the non-minimal 
solution $x_0=70$): 
$x_i=\mid 70+250i\mid $ for $A_i=2+14i+25i^2$. The solutions are $(A_i, x_i; i)$: 
\underline{$(13, 180; -1)$}, $(41, 320; 1), (74, 430;-2), (130, 570; 2)$, etc. 
Now we describe the general method for finding "vertical" relations. 

\vskip 0.2in

{\bf The method of hidden parameters}. We introduce additional parameters (that 
is, the problem must be complicated at first), and then the relations 
between these new parameters must be applied. For the case of $A=a^2+b^2$ 
we suppose, that all four values $l, m, a$ and $b$ are unknown, and
shall seek the solution of condition (6) as  
\begin{equation}
l = gn + r , ~~~ m = dn + h , ~~~ b = pn + z , ~~~ a = qn + u , 
\end{equation}
with new unknown integer parameters $n, g, r, d, h, p, z, q$ and $u$. 
Seemingly, the substitution $n=0$ into final identities is to no avail, 
because it implies only rewriting of symbols, which it is not the case. This 
is the reason, why the notion of "hidden parameters" appears in the title. 

We suppose that $n$ can be arbitrary in an ultimate solution and group 
Eq.(6) in terms of $n$: 
\begin{equation}
~ ~ ~ (n^3): ~ ~ ~ ~ ~ ~ ~ ~ ~ ~ ~ ~ ~ ~ ~ ~ ~ ~ ~ ~ ~  
g^2p - d^2p - 2gdq = 0 , ~ ~ ~ ~ ~ ~ ~ ~ ~ ~ ~ ~ ~ ~ ~ ~ ~ ~ ~ ~ ~ ~ ~ ~ ~ 
\end{equation}
\begin{equation}
~ ~ ~ (n^2): ~~~~~~~ 2grp - 2dhp + g^2z - d^2z - 2ghq - 2drq - 2ugd = 0 , 
\end{equation}
\begin{equation}
~ ~ ~ (n^1): ~~~~~~~ pr^2 - ph^2 + 2grz - 2dhz - 2rhq - 2ugh - 2dru = 0 , 
\end{equation}
\begin{equation}
~ ~ ~ (n^0): ~ ~ ~ ~ ~ ~ ~ ~ ~ ~ ~ ~ ~ ~ ~ ~ ~ ~ ~ ~ 
z(r^2 - h^2) - 2urh = \pm ~1 .~ ~ ~ ~ ~ ~ ~ ~ ~ ~ ~ ~ ~ ~ ~ ~ ~ ~ ~ ~ ~~ 
\end{equation}
Our aim is to find any solution 
of Eq.(6), namely, $a$ and $b$ for distinct parameters $l$ and $m$. Then, using 
"horizontal" relations (13), we shall obtain all solutions $(A_i, x_i)$ 
of Pell's equation. One can put in (15): 
\begin{equation}
p = gdv, ~~~~~ q = {1\over 2}v(g^2 - d^2), 
\end{equation}
where $v$ is an arbitrary unknown integer. Substituting $p$ and $q$ into 
(16), (17) and resolving the linear system of equations in terms 
of $v$ and $u$, one gets: 
\begin{equation}
v = z{g^2 + d^2\over hg^3 + zd^3}~, ~~~ 
u = -{z[hd^3 - rg^3 - 3gd(rd - hg)]\over 2(hg^3 + rd^3)} . 
\end{equation}
It follows from Eq.(18), that $z$ and $u$ are relatively prime integers; 
$z$ is an odd value 
\begin{equation}
z = hg^3 + rd^3 . 
\end{equation}
Substituting $z$ in Eq.(18), one obtains the condition 
\begin{equation}
rd - hg =~ \pm ~ 1 . 
\end{equation}
Thus, the final expressions are as follows: 
\begin{equation}
~ ~ ~ ~~b_0=ngd(g^2+d^2)+hg^3+rd^3,~a_0={n\over 2}(g^4-d^4)-{1\over 
2}[hd^3-rg^3\mp 3gd]
\end{equation}
with a single condition (22). If for arbitrary coprime numbers $r$ and $h$ we 
shall know any $g$ and $d$ from condition (22), then the numbers $a_0$ and $b_0$ can be 
calculated from expression (23); therefore, taking into account the "horizontal" 
relations (13), all solutions of Pell's equation with parameters $l$ and 
$m$ are found. 

The remark on "hidden parameters" is in order. For $n=0$ we have the following 
solution (which can be proved by direct substitution) from expression (23) 
despite the fact, that it seemed to be impossible at a glance (we rename 
$r\rightarrow l, h\rightarrow m$). \\
~\\
{\Large Proposition 3}. {\it The minimal natural solutions $x_i$ of (1) for different} 
$A_i=a_i^2+b_i^2$ {\it can be taken from expressions (5),(10),(11) with }
\begin{equation}
b_0 = mg^3 + ld^3 , ~~~ a_0 = -{1\over 2}[md^3 - lg^3~\mp ~ 3gd] , 
\end{equation}
\begin{equation}
ld - mg =~ \pm ~ 1 . 
\end{equation}

\vskip 0.25in

Selecting the identities, which represent the unity in various forms (and 
comparing with (25)), and using formulae ((5), (10), (11), (24)) derived above, 
the representation of Pell's equation in form of identities can be written, and for 
specific parameters these identities present the solutions which are of interest to us:

(a) Putting $d=1$, we have (see (25)) $l=mg\pm 1$ (for the odd number $m$ the number $g$ 
is odd; for the even $m$ the value $g$ is arbitrary),
$$
b_i = mg^3 + mg \pm 1 + i(2m^2g \pm 2m) ,  
$$
$$
a_i = -{1\over 2}[m - mg^4 \mp g^3 \mp 3g] + 
i(m^2g^2 + 1 \pm 2mg - m^2),
$$
and the three-parametric identity (for arbitrary values of $m,g,i$) is
$$
\Biggl [ (mg^3 + mg \pm 1 + 2m^2gi \pm 2mi)^2 + \Biggr ( {m(g^4 - 1) \pm 
g(g^2 + 3)\over 2} + m^2g^2i \pm 2mgi + i - im^2\Biggr )^2\Biggr ]\times 
$$
$$
[2(m^2 + m^2g^2 + 1 \pm 2mg){2m(mg \pm 1)a_i + (m^2g^2 + 1 \pm 2mg - m^2)b_i}]^2 
+ 1 = 
$$
$$
[2(a_i^2 + b_i^2)(m^2g^2 + 1 \pm 2mg + m^2)^2 - 1]^2 . 
$$
For example, we present (from here) some minimal natural solutions for $A<150$ (see
formulae (5),(10),(11),(24),(25)):\\ 
for $g=3,~ m=2, ~i=-2$ and at the top signs we have $A=89, ~x=53000$;\\
if $g=5,~ m=2,~ i=-6$ (the top signs), then $x=267000$ for $A=73$;\\
for $g=3, m=3, i=-2$ (the bottom signs) we obtain $A=113, ~x=113296$.

For other identities, which can also be written with the help of "a pen and a 
paper" only (with using formulae (5),(10),(11),(24),(25)), we present
parameters $l$ and $m$, $A$ and $x$:

(b) Putting in (25) $d=g-1$ (the bottom sign), it turns out that some identity 
can be written with $l=g(2T+1)+1,~~ 
m=g(2T+1)-2T$; for $T=3,~ g=3$ we have $l=10,~ m=7$, and with 
"a shift" $i=-2$ the minimal natural solution for $A=137$ is $x=519712$.

(c) Let $g=rd\mp 1$ , $m=2dT+1, l=rm\mp 2T$ (to obtain some identity in (25)); 
for $T=3,~ d=2,~ r=2$ (top sign) the parameters are $m=13,~ l=20$, and with 
"the shift" $i=-1$ we find $A=97, ~x=6377352$.

(d) We substitute $g=J+T(J-1), ~d=J-1, ~r=TJ+J+1,~ l=2gn_1+r, ~m=2dn_1+J$ (we 
have an identity from (25)), and 
for $J=-3, ~T=2, ~n_1=-3$ the parameters are $m=21,~ l=58$; the 
latters (with "the shift" $i=13$) describe the case $A=61, ~x=226153980$.

(e) We have $l=n_1+Tm,~ g=1+dT,~ m=n_1d\mp 1$ ((Eq.(25) becomes an identity); it 
follows for $n_1=5,~ d=6, ~T=3$ (the top sign), that $m=29, ~ l=92$; 
the latter parameters with $i=-41$ satisfy the case of $A=149, ~x=2113761020$, 
etc. 

\section{THE CASE OF PRIME $A = 8N + 3$}

{\Large Proposition 4}. {\it Let} $A$ {\it be a prime number} $A=8N+3=a^2+2b^2$; 
{\it then the minimal natural solution of Pell's equation (1) is}
\begin{equation}
x = |~4blm + a|l^2 - 2m^2|~|\cdot (l^2 + 2m^2),
\end{equation}
{\it with} $l$ {\it and} $m$ {\it taken from the distinctive form}
\begin{equation}
b|l^2 - 2m^2| - 2alm = \pm 1 . 
\end{equation}

{\it Proof}. 
The minimal nontrivial solution is an odd number in this case (see the case of 
$A=4N+1$)
\begin{equation}
x = QS.
\end{equation}
Using the simplicity of $A$, two possibilities follow again:
$$
(i) ~~~~~~~~~~~~~~ y - 1 = AS^2 ~, ~~~ y + 1 = Q^2 .
$$
Subtracting the first equation from the second one, it turns
out that
\begin{equation}
Q^2 = AS^2 + 2 .
\end{equation}
However, this equation is contradictory to solve modulo 8.
$$
(ii) ~~~~~~~~~~~~~~~ y + 1 = AS^2 ~, ~~~ y - 1 = Q^2 .
$$
\begin{equation}
AS^2 = Q^2 + 2\cdot 1^2 .
\end{equation}
Since the number in the right-hand side is representable as a sum of
a square and a double square, the number in the left-hand side
has similar representation (the procedure of division). 
We present $A$ as $a^2+2b^2$. Since one of squares is unity, $S$ may be 
represented as $S=l^2+2m^2$ 
and $S^2=(l^2-2m^2)^2+2(2lm)^2$. Using expressions (4),(30), the proposition 
is proved. 
\vskip 0.25in

Note, that this case can be applied to some cases of complex
$A=a^2+2b^2$. Using the expression $y = A(l^2 + 2m^2)^2 - 1$, we see, that
$l,m \le \sqrt[4]{(y+1)/A} < \sqrt[4]{(\sqrt{A}x+2)/A}$. We have considerable decreasing of calculations again (for the method of 
sequential differences), since we use some part of the cycle only and 
calculate comparatively small parameters $l,m$. The question of
representation of $A$ as $a^2+2b^2$ and of choosing the sign (for removing of 
the modulus) can be bypassed with using the inverse calculations (see Section 5 below).

We present the method for decreasing of calculations for the method of 
continued fractions. Using Proposition 4 and the correspondence of the form 
(2) and the expression (9), we obtain the following result. 

(I) If we have 
$$r_j=2r_{j-1} ~ ~ \Rightarrow ~ ~ a=a_j^t, ~ ~ ~ b=r_{j-1}, ~ ~ ~ m=B_{j-2}, ~ ~ ~ l=B_{j-1}. 
$$

(II) If we have 
$$
r_{j-1}=2r_j ~ ~ \Rightarrow ~ ~ a=a_j^t, ~ ~ ~ b=r_j, ~ ~ ~ m=B_{j-1}, ~ ~ ~ l=B_{j-2}. 
$$
Everywhere $\{ B_n\}$ is the recursive sequence given by $B_{-2}=1, B_{-1}=0, 
B_n=d_nB_{n-1}+B_{n-2}$. We demonstrate the case of $A=139$.
$$
\sqrt{139}=11+{18\over \sqrt{139}+11}, ~~ {\sqrt{139}+11\over 18}=1+{5\over 
\sqrt{139}+7}, ~~ {\sqrt{139}+7\over 5}=3+{15\over \sqrt{139}+8},
$$
$$ 
{\sqrt{139}+8\over 15}=1+{\underline{6}\over \sqrt{139}+7}, ~~ 
{\sqrt{139}+7\over 6}=3+{\underline{3}\over \sqrt{139}+11},
$$
We see that (II) $b=3, a=11; B_0=1, B_1=1, B_2=4, l=B_3=5, m=B_4=19$. 
Substituting $l,m,a,b$ into (26), the minimal positive solution $x=6578829$ follows. 
Instead of 19 steps for the continued fraction expansion of $\sqrt{139}$ we use 
5 steps only, and instead of calculation of $B_{17}$ (it is the solution - see 
Introduction), we calculate $B_4$ only. Therefore, the number of calculation 
steps can be considerably decreased.\\
~\\
{\Large Proposition 5}. {\it Let} $x_0$ {\it be the minimal natural solution of Eq.(1) for} 
$A_0=a_0^2+2b_0^2$ {\it and definite parameters} $l, m$ {\it from condition (27), 
then} $x_i$ {\it are the minimal natural solutions of Eq.(1) for} $A_i$ {\it and} $i\ne -1$, {\it where}
\begin{equation}
x_i = |x_0 + iS^3| ,
\end{equation}
\begin{equation}
A_i = A_0 + 2i{x_0\over S} + i^2S^2 .
\end{equation}

{\it Proof}. Designating $k=|l^2-2m^2|,~~~ t=2lm,$
the condition (27) can be rewritten as (6'). 
Because Eq.(6') is satisfied by $a_0$ and $b_0$, all solutions are (13).
The solution can be written as $x_i=|a_ik + 2b_it|S, ~A_i=a_i^2+2b_i^2$.
Because $k^2+2t^2=S^2$, the net formulae follow. 
The proof of solution minimality is valid (see Section 2).
\vskip 0.25in

For example, (using the nonminimal solution $x_0=12$ and parameters $l=1,m=1$) 
we can obtain the minimal natural solutions $x_i=|12+27i|$ for $A_i=2+8i+9i^2$.
The $(A_i,x_i;i)$ solutions (from (31),(32)) are: $(19, 39; 1), (54, 66; 2), 
(107, 93; 3), (22, 42; -2), (59, 69; -3), (114, 96; -4)$, etc.\\
~\\
{\Large Proposition 6}. {\it The solution of Eq.(1) for the case of} 
$A_i=8N+3=a_i^2+2b_i^2$ {\it is 
presented by expressions (26),(31),(32) with}
\begin{equation}
b_0 = 4mg_1^3 + ld^3 , ~~~ a_0 = |md^3 - 2lg_1^3 \mp ~3g_1d| ; 
\end{equation} 
\begin{equation}
ld - 2mg_1 =~ \pm ~1 . 
\end{equation}

{\it Proof}. The direct substitutions prove this proposition. However, we shall 
outline the derivation of the statement with using the method of hidden parameters. 
Substituting (14) into condition (27) (we substitute~ $sign=\pm$ ~ instead 
of $|...|$) and grouping power series in $n$, we have four 
equations again. Substituting 
\begin{equation}
p = sign\cdot gdv , ~~~~ q = {1\over 2}v(g^2 - 2d^2),
\end{equation}
and resolving the linear set of simultaneous equations in terms of $v$ and 
$u$, we find: 
\begin{equation}
v =~ sign\cdot z{g^2 + 2d^2\over hg^3 + 2rd^3}~ , ~~~ 
u =~ -sign\cdot {z[4hd^3 - rg^3 + 6gd(hg - rd)]\over 2(hg^3 + 2rd^3)} . 
\end{equation}
In order that Eq.(27) be solvable, it turns out that 
\begin{equation}
z = {1\over 2}(hg^3 + 2rd^3) , 
\end{equation}
where $g=2g_1$. The final condition is 
\begin{equation}
rd - 2hg_1 =~ \mp sign\cdot 1 
\end{equation}
The solution is: 
\begin{equation}
a_0 =~ sign\cdot [(4g_1^4 - d ^4)n - [hd^3 - 2rg_1^3~ \pm sign\cdot 3g_1d]] , 
\end{equation}
\begin{equation}
b_0 = 2g_1d(2g_1^2 + d^2)n + 4hg_1^3 + rd^3 . 
\end{equation}
"Horizontal" relations (31),(32) provide all solutions for the given 
parameters. Although we have sought again only some of solutions, all 
solutions of Eq.(1) can be obtained with the help of condition (38) 
in this case. Letting again $n=0$ and rename $r\rightarrow l, h\rightarrow m$, 
we have the net result from the proposition.
\vskip 0.25in

For example, some identities for (1) can be obtained in the following way.

(a) If $d=1$, then (from (34)) $l=2g_1m\pm 1$; for $g_1=2,~ m=1$ and the top 
sign one gets $l=5$; for "the shift" $i=-4$ we have (from (26),(31),(32),(33)) 
the case of $A=67, ~x=5967$.

(b) If $g_1=3,~ m=1$ (the top sign), then $l=7$, and with $i=-8$ we have 
the case of $A=118,~ x=28254$.

(c) We substitute in (34): $d=2Tg_1\pm ~1$, and 
$l=2g_1K+1,~ m=Tl\pm ~K$. Let $g_1=2, ~K=1$ 
in the identity; then $l=5$, and for $T=4$ (the bottom sign) 
the parameter $m$ is: $m=19$; with "the shift" $i=-92$ the solution 
follows: $x=6578829$ for $A=139$.

\section{THE CASE OF PRIME $A = 8N + 7$}

On one hand, the solution $x$ cannot be 
an even number (see Section 2). On the other hand, Eq.(30) is intractable to solve 
modulo 8 for $A=8N+7$. Therefore, we obtained the solution (28) with coprime 
odd parameters which can be found 
from the condition (29) by sequential differences. We note, that the 
representation of prime $A=8N+7$ as $a^2-2b^2$ follows from (29) (the 
existence of the solution and the division procedure). Substituting 
$A=a^2-2b^2$ and $S=l^2-2m^2$ into (29); thus the following proposition 
is proved:\\ 
~\\
{\Large Proposition 7}. {\it The minimal natural solution of Eq.(1) for prime 
(and some complex)} $A=a^2-2b^2$ {\it is} 
\begin{equation}
x = |~(a(l^2 + 2m^2) - 4blm)(l^2 - 2m^2)~| , 
\end{equation}
\begin{equation}
2alm - b(l^2 + 2m^2) = \pm 1 .
\end{equation}
\vskip 0.25in
Since $l>m$ (we substitute $m=l-z$, further $l=m+z$), the distinctive form 
(the sum of coefficients at the squares of parameters equals the coefficient 
at the product of parameters) is 
\begin{equation}
(2a - 3b)m^2 - bz^2 + 2(a - b)mz = \pm 1 .
\end{equation}

Writing $y = A(l^2-2m^2)^2 + 1$, we obtain $l,m < \sqrt[4]{y/A} <
\sqrt[4]{(\sqrt{A}x+1)/A}$ (i.e. the parameters $l$ and $m$ are comparatively
small again). Using (2), (9) and Proposition 7, the method for decreasing of calculations 
for the method of continued fractions can be suggested. If we have 
$$
r_{j-1}+
r_j=2a_j^t ~ ~ \Rightarrow ~ ~ b=r_{j-1}, ~ ~ ~ a=a_j^t+r_{j-1}, ~ ~ ~ m=B_{j-2}, ~ ~ ~ l=B_{j-1}+
B_{j-2}
$$
with the recursive sequence $\{ B_n\}$ given by $B_{-2}=1, B_{-1}=0, 
B_n=d_nB_{n-1}+B_{n-2}$. 

One can write the relations between
$(x_0, A_0)$ and $(x_i, A_i)$, so that all substitutions in the method of
sequential differences are identical. Designating
 $t = 2lm, ~ k = l^2 + 2m^2$,
the condition (42) can be rewritten as (6'). 
Because Eq.(6') is satisfied by $a_0$ and $b_0$, all solutions
can be taken from expression (13). 
The solution can be written as $x_i=|~(a_ik - 2b_it)S~| , ~A_i=a_i^2-2b_i^2$.
Because $S^2+2t^2=k^2$, we proved the following proposition:\\ 
~\\
{\Large Proposition 8} {\it All "horizontal" relations between minimal natural 
solutions} $x_0$ {\it for} $A_0$ {\it and} $x_i$ {\it of Pell's equation (1) 
(and} $i\ne -1$ {\it ) for 
the case of} $A_i=a_i^2-2b_i^2$ {\it are presented by the formulae}
\begin{equation}
x_i = |x_0 + iS^3| , ~~~ S = l^2 - 2m^2 , 
\end{equation}
\begin{equation}
A_i = A_0 + 2i{x_0\over S} + i^2S^2 .
\end{equation}
\vskip 0.25in

For example (we use parameters $l=3,m=1$ for nonminimal solution $x=70, A=2$), 
$x_i=|70+343i|, A_i=2+20i+49i^2$, i.e. 
$(A_i,x_i,i):~ (31, 273; -1), (71, 413; 1)$, etc. 

All "vertical" relations can be found by the method of hidden parameters. We 
present the results for expression (14).
$$
u = {4hd^3 - rg^3 \pm 6gd\over 4}, ~~~ z = {2d^3r - hg^3\over 2}, 
$$
$$
q = {g^4 - 4d^4\over 4}, ~~~~ p = {dg(g^2 - 2d^2)\over 2}, ~~~ rd - hg = \pm 1.
$$
For $n=0, ~g=2g_1$ (with $r\rightarrow l, h\rightarrow m$) we have the 
following proposition, which can be ~checked by the direct substitutions.\\ 
~\\
{\Large Proposition 9}. {\it The solution of Pell's equation (1) for} 
$A_i=a_i^2-2b_i^2$ {\it is given by expressions (41),(44),(45) with}
\begin{equation}
a_0 = md^3 - 2lg_1^3 \pm 3g_1d, ~~~~ b_0 = ld^3 - 4mg_1^3, 
\end{equation}
\begin{equation}
ld - 2mg_1 = \pm 1.
\end{equation}
\vskip 0.25in

We present the examples of identities below.

(a) The condition (47) is an identity for $d=1, m=1, l=2g_1+1$. For $g_1=3, 
i=7$ the case $A=103, x=22419$ follows from (41),(44)-(46).

(b) The condition (47) becomes an identity for $d=1, m=2g_1, l=4g_1^2\pm 1$. 
If $g_1=2, i=1$, then (see (41),(44)-(46)) $A=127, x=419775$.

\section{OTHER CASES}

The rest case of $A=p_1p_2$, where $p_1$ and $p_2$
are coprime numbers (however, they can be complex numbers), can be broken 
down into the following ones. 
 
(A) The solution $x$ of (1) is an even number $x=2QS$.
It follows from
$$
y + 1 = 2p_1S^2 ~ , ~~~~ y - 1 = 2p_2Q^2 ,
$$
that the distinctive form is
\begin{equation}
p_1S^2 - p_2Q^2 = 1 .
\end{equation}
For the method of continued fractions, this form does not present in (9). 
However, expression some parameter in terms of other one leads to the form, 
which is presented in (9). For the decreasing of calculations, if we have 
$$
Kr_j=a_j^t, ~ ~ K=1,2,..., ~ ~ \Rightarrow ~ ~ Q=B_{j-1}, ~ ~ ~ S=KB_{j-1}+B_{j-2},
$$ where 
the recursive sequence $\{ B_n\}$ is $B_{-2}=1, B_{-1}=0, B_n=d_nB_{n-1}+
B_{n-2}$.

Since $p_{1}^i = p_{01} + iQ^2 , ~~~ p_{2}^i = p_{02} + iS^2 $, these 
"horizontal" relations give 
\begin{equation}
x_i = 2QS , ~~~ A_i = A_0 + i^2{x_0^2\over 4} + i(p_{01}S^2 + p_{02}Q^2) .
\end{equation}
The following "vertical" relations can be checked by direct substitutions: 
\begin{equation}
p_{01} = m^2(3Ql - Sm) , ~~~ p_{02} = l^2(3Sm - Ql) , 
\end{equation}
with the condition  
\begin{equation}
Ql - Sm = 1 . 
\end{equation}
The solution always exists for coprime numbers $S$ and $Q$; that is, we have 
obtained all relations for $x_i$ and $A^i=p_{01}^ip_{02}^i$ (see (49),(50)). 

(B) The solution of (1) is an odd number $x=SQ$. The condition is
$$
p_1S^2 - p_2Q^2 = \pm 2.
$$
Let $p_1>p_2$, $Q=S+2q$, and one
obtains the distinctive form with equal coefficients at one of squares and 
at the product of parameters
\begin{equation}
{p_1 - p_2\over 2}S^2 - 2p_2q^2 - 2p_2Sq = \pm 1 .
\end{equation}
For the decreasing of calculations in the method of continued fractions, we 
seek the following condition in (9): 
$$
r_j=2a_j^t ~ ~ \Rightarrow ~ ~ x=B_{j-1}(
B_{j-1}+2B_{j-2})
$$
with the well known recursive sequence $\{ B_n\}$.

The "horizontal" relations 
\begin{equation}
x_i = SQ , ~~~ A_i = A_0 + i(p_{01}S^2 + p_{02}Q^2) + i^2x_0^2 
\end{equation}
with the "vertical" relations
\begin{equation}
p_{01} = 2m^2(3Ql - mS) , ~~~ p_{02} = 2l^2(3mS - Ql), 
\end{equation}
\begin{equation}
lQ - mS = 1
\end{equation}
present all solutions of Pell's 
equation for this case (and can be proved by direct substitutions). 

We write the binary quadratic form (2) with $b_n>0, c_n>0$ and 
\begin{equation}
b_n - c_n + 2a_n = 1~
\end{equation}
(the sequential differences lead to $X_n=Y_n=1$ for some step $n$).
The net formulae are
\begin{equation}
a_n = A - k(k+1) .
\end{equation}
\begin{equation}
b_n = (k+1)^2 - A ,
\end{equation}
\begin{equation}
c_n = A - k^2 .
\end{equation}
It follows from $b_n>0 , c_n>0$, that
\begin{equation}
k^2 < A < (k+1)^2.
\end{equation}
This condition for the given $A$ (not a square)
uniquely determines the value $k$ and, therefore, the values
$a_n , b_n , c_n$. It can be proved that (2) can be reduced to (1). 

\vskip 0.2in

{\bf Inverse calculations} (for the method of sequential differences). 
We believe that the "coordinates" of a given 
number $A$ are known (the nearest square $k^2$ and the "distance" from it).  
Therefore, using (57)-(60), one can write the "ultimate form", which
represents the "one" (by substitutions $X_n=1, Y_n=1$). 
Using the symmetry of cycle with respect to the forms with $+a_i$ and $-a_i$, 
we change the sign of $a_n$ and obtain the "start form". By the sequential 
differences method we obtain the "final form". This is one of
five types of the form discussed above:

(I) the coefficients at the squares of parameters are equal (see (6));

(II) the coefficients at the squares differ twice from each other (see (27));

(III) the sum of coefficients at the squares of parameters equals the 
coefficient at the product of parameters (see (43));

(IV) the product of parameters is missing from the form (see (48));

(V) the coefficient at one
of squares equals the coefficient at the product of
parameters (see (52)).\\
Changing the sign of the product of parameters, we obtain the "distinctive
form" (with real sign), that is, the real sign in the solution follows. To 
find the "distinctive parameters", it is necessary:

1. to find the "start parameters" $X_n$ and $Y_n$ in terms of
the "final parameters";

2. to express the "final parameters" in terms of the "start parameters";

3. to change the sign at one of the "start parameters"
(either $X_n$ or $Y_n$) in both expressions, to substitute
$X_n=Y_n=1$ and to take the magnitudes. As a result, we obtain
the "distinctive parameters" (and the solution of (1) for the given 
number $A$).
Notice, that beginning from the "start form" we can also found the 
representation of number $A$ as $a^2+b^2$, or as $a^2+2b^2$, or the expansion 
of $A$ into two coprime factors (the same result can be obtained in continued 
fraction expansion).

For example, we outline the case of $A=103$. The "start form"
($10^2<103<11^2,k=10$) is 
$18Y_n^2-3X_n^2+14X_nY_n=1$. We have the "final form" (III) at the 5th step: 
$13m_1^2-3z_1^2-16m_1z_1=1$. We obtain $b=3, ~a=11, ~ X_n=m_1,~ Y_n=5X_n+z_1$. Therefore, 
$z=6, m=1, l=7$ and it follows from (41): $x=22419$.

\section{SOME RELATED PROBLEMS}

The methods described above permit to explicitly 
solve some related problems in the general case (for example, $y^2-Ax^2=\pm j^2;
\pm 2j^2$, etc.). We consider the case (see [7]):
\begin{equation}
y^2 - Ax^2 = -3.
\end{equation}
Since there exists the procedure of division for the numbers $y^2+3j^2$, two 
possibilities follow for odd $A$ (the case of $4A$ can also be 
reduced to them).

(i) The minimal natural solution $x$ of (61) is an odd number for $A=a^2+3b^2=4N+3$ 
(with an odd number $b$) in the first case: 
\begin{equation}
x = l^2 + 3m^2, ~~~  y = |~(a|l^2 - 3m^2| - 6blm)~| , 
\end{equation}
with the condition 
\begin{equation}
2alm - b|l^2 - 3m^2| = \pm 1.
\end{equation}
To decrease calculations for the method of continued fractions, we obtain the 
following result. 

(I) If we have 
$$
r_j=3r_{j-1} ~ ~ \Rightarrow ~ ~ l=B_{j-1}, 
~ ~ ~ m=B_{j-2}, ~ ~ ~ a=a_j^t, ~ ~ ~ b=r_{j-1};
$$

(II) if we have 
$$
r_{j-1}=3r_{j} ~ ~ \Rightarrow ~ ~ l=B_{j-2}, ~ ~ ~ m=B_{j-1}, ~ ~ ~ b=r_j, ~ ~ ~ a=a_j^t
$$
with the recursive sequence $\{ B_n\}: 
B_{-2}=1, B_{-1}=0, B_n=d_nB_{n-1}+B_{n-2}$. 

All "horizontal" relations can be taken from
\begin{equation}
a_i = a_0 + i|l^2 - 3m^2| , ~~~~ b_i = b_0 + 2ilm, ~~~ A_i = a_i^2 + 3b_i^2. 
\end{equation}
Using the method of hidden parameters, one can find the parametric solution 
of this problem (which can be proved by direct substitutions): 
\begin{equation}
a_0 = {3\over 2}[d^3m + 9dmg^2 - 3d^2lg - 3lg^3], ~~~~~~ b_0 = 9g^3m + 
ld^3, 
\end{equation}
where the condition for $d$ and $g$ is 
\begin{equation}
dl - 3gm = \pm 1 . 
\end{equation}

(ii) The second case with an even solution of (61): $x=2(l^2+3m^2)$, takes place for 
$A=A_0^2+3B_0^2=4N+1$ ($B_0$ is an even number). In this case 
($x\rightarrow 2x_1,~ A\rightarrow 4A_1$) all previous formulae 
are valid with substitutions 
\begin{equation}
a = A_0\pm 3B_0, ~~~~ b = A_0 \mp B_0. 
\end{equation}

We consider the case $A=1729=4N+1$~ (see [7], for comparison). 
$$
\sqrt{4\cdot 1729}=83+{27\over \sqrt{6916}+83}, ~~ {\sqrt{6916}+83\over 27}=
6+{\underline{25}\over \sqrt{6916}+79}, 
$$
$$ 
{\sqrt{6916}+79\over 25}=6+{\underline{75}\over \sqrt{6916}+71}, 
$$
We see that (I) $m=B_1=6, l=B_2=37, x=2954, y=122831$. We outline the 
parametric identity which includes this case (see (62),(64),(65)): the 
condition (66) is an identity for $d=1, l=3gm+1$ (for $m=6, g=2, i=-1$ 
one obtains from (64),(65): $b=25, a=71$, and the case follows from (67)). 

\section{CONCLUSIONS}

The general case of Pell's equation (1) can be broken down into the special 
cases, and in each of these cases the number of calculation steps can be
considerably lowered for the method of sequential differences and for the 
method of continued fractions (the minimal natural solution can be found with the 
help of "a pen and a paper" within a few minutes). The item-by-item examination
of particular cases may seem increasing the volume of
calculations. However, it is not the case. The author computed Pell's equation
(in Python): (a) with using the method of continued fraction
(standard algorithm), and (b) with using (in addition) the particular cases
and formulas from this article. We have the following results. For all
$A<1000000$ it is necessary about 3 hours and 2 minutes for the standard algorithm and
1~h. 43~m. for the modified algorithm; for all $A<2000000$ it is required about
8~h. 50~m. for the
standard algorithm and 4~h. 50~m. for the modified algorithm.
Thus, the lowering of computation time
is not less than 1.828 times in the last case (and 1.827 times, if we seek and write the absolute
maxima also). Since the maxima significantly increase with increasing of $A$,
this result can be
improved (the basic time surplus we have from the most difficult cases).
Computer programs (from Internet) use the method of continued fraction. The
most fast algorithm (in Mathematica) uses additionally one particular case
only: $y_1^2 - Ax_1^2 = -1$ (factually, Eq.(8) from Section 2, but not Eq.(6)).
Therefore, the suggested algorithm is more effective in this particular case
also. The methods, which 
permit to decrease the scope of calculations for finding the table of 
minimal natural solutions of Pell's equation (1), are investigated in the work. Using 
these methods, it is possible to establish the relations (to deduce the
appropriate formulae) between a solution $x_0$ for $A_0$ and the infinite 
number of minimal natural solutions $x_i$ for $A_i$ in various cases.

The establishment of "horizontal" (both basic parameters are fixed) and "vertical" 
(basic parameters are variable) relations permits to write a 
parametric representation of Pell's equation (1) as identities. With the help of 
"the method of hidden parameters" all relations can be found. This 
method consists in the following procedure: 
first, it is necessary to complicate the problem by inserting an additional 
number of new parameters; then, by obtaining relations between these new 
parameters, the appropriate solutions can be found; and, finally, by reducing the
problem to the initial number of parameters, the solutions can be represented 
as identities by means of initial parameters. Here these relations permit 
to find the minimal natural solution of Pell's equation (1) (for the given $A_0$) and to 
present Pell's equation (1) as a parametric identity in the explicit form (for
infinite natural numbers of $(A_i, x_i)$, including $(A_0, x_0)$). 
These
relations do not connected with composition in cycles, but Pell's equation (1)
is connected with indeterminate linear equation $ax - by = 1$.

{}
\end{document}